\documentclass[12pt]{article}
\usepackage{latexsym}
\begin{document}
\newcounter{thm}
\newtheorem{Def}[thm]{Definition}
\newtheorem{Thm}[thm]{Theorem}
\newtheorem{Lem}[thm]{Lemma}
\newtheorem{Cor}[thm]{Corollary}
\newtheorem{Prop}[thm]{Proposition}
\newcommand{\vlimsup}{\mathop{\overline{\lim}}}
\newcommand{\vliminf}{\mathop{\underline{\lim}}}
\newcommand{\Av}{\mathop{\mbox{Av}}}
\newcommand{\spec}{{\rm spec}}
\newcommand{\qed}{\hfill $\Box$ \\}
\newcommand{\subarray}[2]{\stackrel{\scriptstyle #1}{#2}}
\setlength{\baselineskip}{15pt}
\def\textmc{\rm}
\def\({(\!(}
\def\){)\!)}
\def\R{{\bf R}}
\def\Z{{\bf Z}}
\def\N{{\bf N}}
\def\C{{\bf C}}
\def\T{{\bf T}}
\def\E{{\bf E}}
\def\H{{\bf H}}
\def\Prob{{\bf P}}
\def\M{{\cal M}}     
\def\F{{\cal F}}
\def\G{{\cal G}}
\def\D{{\cal D}}
\def\X{{\cal X}}
\def\A{{\cal A}}
\def\B{{\cal B}}
\def\L{{\cal L}}
\def\a{\alpha}
\def\b{\beta}
\def\e{\varepsilon}
\def\de{\delta}
\def\ga{\gamma}
\def\k{\kappa}
\def\la{\lambda}
\def\fa{\varphi}
\def\th{\theta}
\def\si{\sigma}
\def\t{\tau}
\def\om{\omega}
\def\De{\Delta}
\def\Ga{\Gamma}
\def\La{\Lambda}
\def\Om{\Omega}
\def\Th{\Theta}
\def\lan{\langle}
\def\ran{\rangle}
\def\lbr{\left(}
\def\rbr{\right)}
\def\const{\;\operatorname{const}} 
\def\dist{\operatorname{dist}} 
\def\Tr{\operatorname{Tr}}
\def\quadd{\qquad\qquad}
\def\n{\noindent}
\def\beq{\begin{eqnarray*}}
\def\eeq{\end{eqnarray*}}
\def\supp{\mbox{supp}}
\def\beqn{\begin{equation}}
\def\eeqn{\end{equation}}
\def\bp{{\bf p}}
\def\br{{\bf r}}
\def\v2{\vskip2mm}
\def\pf{{\it Proof.}}

\begin{center}
{\Large The hitting distribution of a line segment   \\ for   two dimensional random walks} \\
\vskip6mm
{K\^ohei UCHIYAMA} \\
\vskip2mm
{Department of Mathematics, Tokyo Institute of Technology} \\
{Oh-okayama, Meguro Tokyo 152-8551\\
e-mail: \,uchiyama@math.titech.ac.jp}
\end{center}

{\bf Abstract}\,\,  Asymptotic estimates of the hitting distribution of a long segment on the real axis for two dimensional random walks on $\Z^2$ of zero mean and finite variances are obtained: some are general and exhibit its apparent similarity to  the corresponding Brownian density, while  others are so detailed as  to involve  certain characteristics of the random walk.  \footnote{{\it key words}:  harmonic measure in a  slit plane, line segment,   asymptotic formula, random walk of zero mean and finite variances. }

\section{Introduction and Results}

Let $S_n = z+ \xi_1+\cdots+ \xi_n$ be a two dimensional random walk  of i.i.d. increments  $\xi_1, \xi_2, \ldots$ and initial position $S_0=z$ moving on the square lattice $\Z^2$, which we suppose to be  embedded  in the complex plane $\C$.  
 Let $n$ be a positive integer and  denote by $H^{I(n)}_z(s)$  the probability that the first visit  (after time $0$) to the interval 
$\{-n+1,\ldots,n-1\}$ of   the random walk  $S_{\cdot}$ starting
    at $z$ takes place at $s$.  For the later use it is convenient to define a positive number $n_*$ and an interval $I(n)$ by
$$I(n)=(-n_*,n_*)=\{u\in\R: |u|<n_*\},  \quad n_*=n-1/2.$$
Then $H^{I(n)}_z(s)$, $s\in I(n)$, is written as
$$H^{I(n)}_z(s)=P_{z}[\,\exists j\geq 1,\, S_j=s \,\,\mbox{and}\,\, S_k\notin I(n)\,\mbox{for}\,1\leq k<j\,],$$
where $P_{z}$ stands for the probability of the walk starting at $z\in \Z+i\Z$. An explicit expression  of the corresponding distribution for Brownian motion
is readily derived from the Poisson kernel for the unit disc  in view of the conformal invariance of harmonic measures.  Let $h_x^{I(n)}$  denote  the Brownian analogue of $H^{I(n)}_z$, namely the density of hitting distribution of the interval
$I(n)$ for the two dimensional standard Brownian motion starting at $z$. Then, for $x\in \R\setminus [-n_*, n_*]$
 \beqn\label{000}
 h^{I(n)}_x(s)= \frac{\sqrt{x^2-n_*^2}}{\pi|x-s|}\cdot\frac1{\sqrt{n_*^2-s^2}} \quad \quad (s\in I(n))
 \eeqn
 (see Appendix (A), in  which we compute   $h^{I(n)}_z(s)$ for general  $z$). From the Donsker's invariance principle it is expected that  $H^{I(n)}_z(s)$  behaves similarly to 
 $ h^{I(n)}_z(s)$ if the covariance matrix of $\xi_1$  is  isotropic, but it is not clear at all  in what sense  they are similar. 
   In the present paper  we compute exact asymptotic  forms of $H^{I(n)}_x(s)$  for all $x\in \Z$ as $|x-s|\wedge n \to\infty$ : the first  of them  exhibits its apparent similarity to  $ h^{I(n)}_x(s)$ and  the others give finer  estimates that  involve  certain characteristics of the random walk. 
  The case of  non-real  initial  sites  will be    briefly discussed   in Appendix (D).  
  
  The problem with real initial sites  reduces  to a one-dimensional one.   Indeed, if $S$ starts at a point $x\in \Z\subset \C$ and  $X= (X_n)$ is  its trace  on the real axis, namely $X$ is a  one-dimensional random walk on $\Z$ imbedded in $(S_n)$ with $X_n$ being  the position of  the $n$-th return  of $S$ to the real axis, then  $H_x^{(n)}(\cdot)$ equals  the hitting distribution of $I(n)$ for $X_n$.  It is remarked that   the increment distribution of $X$ is almost Cauchy in the sense that  its tails are asymptotically  $C/|x|$ in both directions \cite{U1}. 
     
 For the symmetric simple random walk H. Kesten has obtained   the upper bound 
\beqn\label{K}
H^{I(n)}_\infty(s):=\lim_{|z|\to\infty} H^{I(n)}_z(s)\leq C[n(n-s)]^{-1/2} \quad\quad  (0\leq s <n)
\eeqn
in \cite{K} (the limit on the left-hand  side  exists  (\cite{S}:Theorem 14.1, p.141)) and applied it to a study of the DLA model in \cite{K1} (cf. also \cite{K2}; a unified exposition is found in \cite{L}).  For a rectangle with a side on the real axis Lawler and Limic \cite{LL} give  an explicit expression for  the hitting distribution of its boundary  for a  simple random walk started inside it and, by taking limits, derive from it the corresponding ones for a half-infinite strip and a quadrant. For a quadrant of the plane, one half of it split  along its diagonal line and the complements of these regions as well  Fukai \cite{F} obtains very detailed evaluations of the hitting distributions by exploiting the properties special to  simple random walk.

Throughout this paper we suppose that the walk $S_n$ is irreducible, $E_0[S_1]=0$ and
\beqn\label{de}
E_0[|S_1|^{2+\de}] <\infty \quad \mbox{either for\,\, $\de=0$\,\, or for some} \quad\de>1/2;
\eeqn
we  make an  explicit  reference  to $\de$ in the latter case, while no reference is  understood to mean   the case $\de=0$. 
\begin{Thm}\label{thm1} \,\, Let $\de>1/2$ in  (\ref{de}).  Then uniformly for integers $s\in I(n)$ and $x, |x| \geq n$,  as $n\to\infty$ 
\beqn\label{(1)}
H^{I(n)}_x(s)= h^{I(n)}_x(s)\Bigg[1+O\bigg(\frac1{\sqrt{(|x| -n_*)\wedge(n-|s|)}}\bigg)\Bigg].
\eeqn
\end{Thm}
\v2

 From  Theorem \ref{thm1} it  follows that 
 \beqn\label{infty}
 H^{I(n)}_\infty(s) =\frac{1}{\pi }\cdot\frac1{\sqrt{n_*^2-s^2}}\Bigg[1+O\bigg(\frac1{\sqrt{n-|s|}}\bigg)\Bigg] \quad \quad (s\in I(n));
 \eeqn 
  (\ref{K}) is thus  refined.  Indeed,   
the probability that the walk starting at $z$ hits the real line in the interval  $I(N)$ tends to zero  for any $N>n$ as $|z|\to \infty$ so that  $H^{I(n)}_\infty(s)$  is represented as the limit of a convex combination of  $H^{I(n)}_x(s)$, $|x|> N$,  which with, e.g.,  $N=2n$ shows the relation above in view of (\ref{(1)}).

When either $|x|-n$ or  $n-|s|$  remains bounded,  (\ref{(1)}) does not determine the asymptotic form of $H^{I(n)}$.  The next theorem improves the estimate   in this respect in the case $\de=0$; in particular it determines  a precise asymptotic form of $H^{I(n)}_\infty(s)$ valid uniformly for  $s\in I(n)$, which  is not provided by (\ref{infty}).  (See
Section 4 (Theorems \ref{thm5}, \ref{thmII-2}) for the case $\de>1/2$.)  The result is expressed by means of  a pair of renewal functions, $\mu(y)$ and $\nu(y)$, $y\geq 0$, associated with the imbedded random walk  $(X_n)$ --- the trace of $S$---  on the real line mentioned above. They are characterized as positive solutions of the Wiener-Hopf equations
$$\mu(y)= E_{-y}[\mu(-X_1); X_1\leq 0]\quad\mbox{and}\quad \nu(y)= E_{y}[\nu(X_1); X_1\geq 0]$$
together with the pairing condition    $\mu(0)\nu(0) =\pi \sigma^{-2}e^{ \sum_{k=1}^\infty k^{-1}P_0[ X_k=0]}$, except for determination of   $\mu(0)>0$ (or $\nu(0)$) that is in our disposal:  we are to  single out $\mu(0)$ appropriately for the present purpose.  Here $\sigma^2$ is the square root of the determinant of the covariance matrix $Q$  of the  i.i.d. increments $\xi_k$: $\sigma := (\det Q)^{1/4}$:  the   quadratic form of $Q$  is  given   by  $E_0[(S_1\cdot \th)^2]$. The equations above plainly  say that $\mu$ and $\nu$  are harmonic for, respectively,  the walks  $-X$ and $X$  killed on hitting the negative  half-line.   We extend $\mu$ and $\nu$ to $y<0$ by these equations.
 It then follows that $\mu(y)$ and $\nu(y), \,y\in \Z$  are (strictly) increasing. 
  We can and do choose  $\mu(0)$ so that
\beqn\label{mu0}
  \frac{\mu(y)}{\sqrt{y}}\,\longrightarrow\,\frac{2}{\sigma^2}   \quad \mbox{and} \quad
 \mu(-y)\sqrt y \,\longrightarrow \,1\quad \mbox{as}\quad  y\to\infty,
\eeqn
which entails the same property for $\nu$ in place of $\mu$ (\cite{U}: Theorem 1.1). 
 (For more details  see Appendix (C).) 

\begin{Thm}\label{thm2} 
\, {\rm (i)}\,   Uniformly for  $0\leq s<n$ and $x\geq n$, as $n\to\infty$ and  $x-s\to \infty$  
\[
H^{I(n)}_x(s)=\frac{\sigma^2}{2\pi}\cdot\frac{\nu(x-n)\mu(-n+s)}{x- s}\cdot\sqrt{\frac{x+n}{n+s}} \Big(1+o(1)\Big).
\]
{\rm (ii)}\,\, Uniformly for  $-n<s\leq 0$ and $x\geq n$, as  $n\to \infty$  
 \[
H^{I(n)}_x(s)=\frac{\sigma^2}{2\pi}\cdot\frac{\nu(x-n)\nu(-n-s)}{x-s}\cdot\sqrt{\frac{x+n}{n-s}} \Big(1+o(1)\Big).
\]
\end{Thm}
\v2

Theorem  \ref{thm2} describes  the asymptotic behavior of  $H^{I(n)}_x(s)$ in the case $x\geq n$; by symmetry we have a similar  result in the case    $x\leq  -n$, actually 
 a translation of  Theorem \ref{thm2} in view of  the  duality of $\mu$ and $\nu$.
\v2\n
{\bf Theorem $2'$}\, {\rm (i)}\,  {\it Uniformly for  $-n < s\leq 0$ and $x\leq -n$, as $n\to\infty$ and  $s-x  \to \infty$  
\[
H^{I(n)}_x(s)=\frac{\sigma^2}{2\pi}\cdot\frac{\mu(-x-n)\nu(-n-s)}{ s-x}\cdot\sqrt{\frac{-x+n}{n-s}} \Big(1+o(1)\Big).
\]
{\rm (ii)}\,\, Uniformly for  $0\leq s <n$ and $x\leq -n$, as  $n\to \infty$ } 
 \[
H^{I(n)}_x(s)=\frac{\sigma^2}{2\pi}\cdot\frac{\mu(-x-n)\mu(-n+s)}{s-x}\cdot\sqrt{\frac{-x+n}{n+s}} \Big(1+o(1)\Big).
\]
\v2
By using the asymptotic form of  the hitting 
distribution of  the real line we can readily deduce  asymptotic forms  of   $H_z^{I(n)}(s)$ for $z\notin \R$ from Theorems  \ref{thm2} and {\bf$2'$} (see Appendix (D)). Here we record only the case when $z=\infty$.
 Since $1/\sqrt{n+s}$ may be replaced by $\mu(-n+s)$ as $n\to\infty$ we obtain the following 
 
\begin{Cor}\label{cor1} \,\, Uniformly for   $s\in I(n)$, as $n\to\infty$  
\[
H^{I(n)}_\infty(s)=\pi^{-1}\mu(-n+s)\nu(-n-s)(1+o(1)).
\]\end{Cor}  

 From Theorems \ref{thm2} and $2'$ we obtain  the second corollary: 
 
 \begin{Cor}\label{cor2} \,\, Uniformly for integers $n\geq 1$, $s\in I(n)$ and $x\in \Z\setminus I(n)$, $H^{I(n)}_x(s)\asymp h_x^{I(n)}(s),$ namely there exists a positive constant
  $C$ independent of $n, s$ and $x$ such that 
  $$C^{-1}h_x^{I(n)}(s)\leq H^{I(n)}_x(s)\leq Ch_x^{I(n)}(s).$$ 
 \end{Cor}

The next theorem provides the asymptotic form of $H^{I(n)}_x(s)$ when $x\in I(n)$.  In view of the corresponding result for the first  visit of the real axis, that may reads
$$P_{z}[\,\exists j\geq 1,\, S_j=s \,\,\mbox{and}\,\, S_k\notin \R\,\mbox{for}\,1\leq k<j\,] \sim  \sigma^2 \lim_{y\to 0}\frac{1}{|y|} h_{x+iy}(s)$$
with  $h_z(s)= |y|/\pi (y^2+(x-s)^2)$ (see  (\ref{5.00})),
 we extend $h^{I(n)}_x(s)$ to the variables $x\in I(n)$, $x\neq s$
by 
$$
 \quad h^{I(n)}_x(s) =\lim_{y\to 0} \, \frac1{|y|}h^{I(n)}_{x+iy}(s).
$$
In Appendix (B) we compute this limit and  find  that 
$$ \quad h^{I(n)}_x(s) = \frac{n_*^2-xs}{\pi (x-s)^2 \sqrt{(n_*^2-x^2)(n_*^2-s^2)}} \quad \quad x,s \in I(n), \, x\neq s.
$$
 (See (\ref{A4}); 
also   the identities (\ref{A3}) and    (\ref{5.0}) for  an underlying idea.)  Let  $S_1^{(1)}$ and  $S_1^{(2)}$ be the real   and imaginary parts of $S_1$, respectively and  let $Y$ be the component of $S^{(1)}$ that is perpendicular to 
$S_1^{(2)}$ under $P_0$, namely $Y= S_1^{(1)}-\om S_1^{(2)}$ where  $\om = E_0[S_1^{(1)}S_1^{(2)}]/E_0[(S_1^{(2)})^2]$.
\begin{Thm}\label{thm21} \,\, Let $Y$ be as above and suppose  the moment  condition
\beqn\label{log_mom} 
E_0\Big[|Y|^2\log |Y|\Big]<\infty. 
\eeqn
 Let  $x,s \in I(n)$. Then 
\v2\n
{\rm (i)} \, as   $(n-|s|)\wedge (n-|x|) \wedge |x-s| \to \infty$  
\[
H^{I(n)}_x(s)=\sigma^2 h^{I(n)}_x(s)(1+o(1));
\]
\v2\n
{\rm (ii)}\,  if $s<x$,  as $(n-x)/(n-s)\to 0$
 \[
H^{I(n)}_x(s)=\frac{ \sigma^2}{ \pi}\cdot \frac{\nu(-n+x)   \nu(-n-s) \sqrt n }{ \sqrt2\, (x-s)^{3/2}} \Big(1+o(1)\Big);
\]
\v2\n
{\rm (ii$'$)}\, if $s>x$,  as $(n-s)/(n-x)\to 0$
 \[
H^{I(n)}_x(s)=\frac{ \sigma^2}{\pi}\cdot \frac{ \mu(-n+s) \mu(-n-x) \sqrt n\, }{\sqrt2\, (s-x)^{3/2}} \Big(1+o(1)\Big).
\]
\end{Thm}

\v2
Observe first  that   the condition $(n-x)/(n-s)\to 0$ in  (ii)  entails 
$$x-s\to \infty, \quad\frac{x}{n} \to 1,\quad  \frac{x-s}{n-s}\to 1  \quad \mbox{ and}\quad  \frac{n^2-xs}{n(n-s)}\to1,$$  and then that     the formula of (ii) implies (and is actually finer than) the formula of (i).  Under the condition $|x-s|\to\infty$ the cases (ii) and (ii$'$) together exhaust  the case when   $(x\vee s)/n\to 1$.  We have an obvious analogue  for (ii$'$), which is a  dual statement of (ii).  Also observe that $h^{I(n)}_x(s)$ as well as  $H^{I(n)}_x(s)$ is bounded away from zero and infinity whenever $|x-s|$  is bounded above by any constant.  These observations  lead  to the following corollary. 

 \begin{Cor}\label{cor3} \, Suppose $E_0[\,|Y|^2\log |Y|\,]<\infty$.  Then, for $x,s \in I(n)$, 
 $C^{-1}h^{I(n)}_x(s) \leq H^{I(n)}_x(s) \leq Ch^{I(n)}_x(s).$
 In particular  if $x$ is kept within  any bounded distance from $n$ and so is $s$  from $-n$,  then  
 $H^{I(n)}_x(s) \asymp 1/n$.
 \end{Cor}

{\sc Remark 1.}  Under the same supposition as in Theorem \ref{thm21} the formulae obtained above can be extended to the general  starting positions $x+iy$ as in \cite{U} but with the resulting formula somewhat complicated (see (\ref{A1})). 

\v2
Denote by $H^{+}_{z}(s)$ the probability that the first visit (after time $0$) to the positive real axis of the walk starting at $z\in \C$ takes place at $s\in \{1,2,3,\ldots\}$:
$$H^{+}_{z}( s)=P_{z}[\,\exists j\geq 1,\, S_j=s \quad \mbox{and}\quad  S_k\notin \{1,2,3,\ldots\} \,\,\mbox{for}\, \,1\leq k<j\,].$$
 Similarly let
 $H^-_{z}(s)$  denote the distribution of the first visiting sites (after time $0$) of the set  $\{-1,-2,-3,\ldots\}$.
The proofs of Theorems \ref{thm1} and \ref{thm2}  rest on the results on $H_x^{\pm}(s)$ obtained in \cite{U} (Theorem 1.1; see also \cite{Uer} for (\ref{(004)}))  as given in the following  theorem (and  also in  (\ref{thmA}), (\ref{thmA0}) and (\ref{thmA1}) later). 

\v2\n
{\bf Theorem} (\cite{U}, \cite{Uer}) \quad  {\it Let $s<0$. Then for $x\geq 0$,  as  $x\vee (-s)\to \infty$  
  \beqn\label{(003)}
 H_x^-(s)=\frac{\sigma^2}{2\pi}\cdot\frac{\nu(x)\mu(s)}{|x-s|} \,\Big(1+o(1)\Big). 
\eeqn
  If  $E_0[|Y|^2\log |Y|] <\infty$  in addition,   then  as $|x-s|\to \infty$ under  $x<0, s<0$, }
  \beqn\label{(004)}
 H_x^-(s)=\frac{\sigma^2}{2\pi}\cdot\frac{|x+s|\nu(x)\mu(s)}{|x-s|^2} \,\Big(1+o(1)\Big). 
\eeqn
\v2
It is warned that in \cite{U} the condition $E_0[|S_1^{(1)}|^2\log |S^{(1)}_1|]<\infty$ is errorneously  assumed where it should be  (\ref{log_mom}) as in Theorem above.
\v2

Comparing  the  formulae given in Theorem with those in Theorems \ref{thm2}(i) and  \ref{thm21}  we find the quite reasonable conclusion  that  $H^{I(n)}_x(s)/ H_{x-n}^-(s-n)\to 1$ if and only if   $x/n\to 1$ and $s/n\to1$ (whether $x$ is larger than $n$  or  not). 

For the symmetric simple random walk (i.e., $P_0[S_1=x]=1/4$ for $x\in\{\pm1, \pm i\}$) we can improve the error estimate in \cite{U} for $H_x^{\pm}(s)$  and accordingly that  in Theorem \ref{thm1} for $H^{I(n)}_x(s)$. 

\begin{Prop}\label{prop0}\,   Let $S_n$ be the symmetric simple random walk. Then
\beqn\label{B-1}
H_x^-(s)=\frac1{\pi}\cdot\frac{1}{x-s}\sqrt {\frac{x\vee1}{-s}}\times\Bigg[1+O\bigg(\frac{1}{{-s}}\bigg)+O\bigg(\frac1{{x\vee1}}\bigg)\Bigg] \quad  (x\geq 0,\,s<0) 
\eeqn 
and as $n\to\infty$
\beqn\label{B-2}
H^{I(n)}_x(s)= h^{I(n)}_x(s)\Bigg[1+O\bigg(\frac1{(|x|-n_*)\wedge(n-|s|)}\bigg)\Bigg]  \quad (|x|\geq n,\,s\in I(n)).
\eeqn
Here the error estimates are uniform for integers $x, s$  subject to the respective constraints indicated in  parentheses.
\end{Prop}

 Theorem \ref{thm1} is proved in Section 2 by taking  for granted certain several results whose proofs are given in Section 3. Proof of Proposition \ref{prop0} is given at the end of Section 3. In Section 4 we make detailed estimation of $H^{I(n)}_x(s)$ when $x$ or $s$ are near the edges of $I(n)$ under the moment condition (\ref{de}); in particular Theorem 2 is proved. Theorem \ref{thm21} is proved in Section 5.

 
\section{Proof of Theorem \ref{thm1}}
The proof of  Theorem \ref{thm1}  primarily rests on the asymptotic  estimates obtained   in \cite{U} of hitting distributions for  half-real-lines.  The first visit of $I(n)$ occurs after the  consecutive  overshoots in the first $k$ alternating  entrances of the walk $X$ into the half-lines $(-\infty, -n)$ and $(n, \infty)$ for some $k=0, 1, 2,\ldots$ and we simply sums up the relevant probabilities over $k$.  Motivated by this we bring in  the sequence of probability kernels in  below. 

Throughout this section we pick  and fix a (large) positive integer $n$, which we shall not designate  in  the notation  introduced in this section even though it depends on  $n$. Let ${\bf 1}(S)$ stand for  the indicator of a statement $S$: ${\bf 1}(S)=1$ or $0$ according as $S$ is true or not. 
Define for integers  $x\geq n$ and $y>-n$,
$$Q(x,y)=\sum_{s=-\infty}^{-n}H^-_{x-n}(s-n)H^+_{s+n}(y+n),$$
$$Q_{I}(x,y)=Q(x,y){\bf 1}(-n<y < n),$$
$$ K_I(x,y)=H^-_{x-n}(y-n){\bf 1}(-n<y <n),$$
and $Q^0={\bf 1}$ (the identity matrix), $Q^1=Q$ and  inductively
\beqn\label{QQ}
Q^k(x,y)=\sum_{u=n}^\infty Q^{k-1}(x,u)Q(u,y) \quad\quad (k=1,2,\ldots),
\eeqn
and finally
$$\La(x,y)=\sum_{k=1}^ \infty Q^k(x,y){\bf 1}(y\geq n).$$
Then
for $x\geq n,$ $ -n<s<n$,
\begin{eqnarray}
\label{QQQ}
H^{I(n)}_x(s)&=&(1+\La) (Q_I+K_I)(x,s) \nonumber\\
&=&\sum_{y=n}^\infty [{\bf 1}(y=x)+ \La(x,y)][Q_I(y,s)+K_I(y,s)].
\end{eqnarray}
 
  The kernels  $Q, K_I, Q_I$ and  $\La$  are probabilities with self-evident meaning. We are to compare them with the corresponding ones, denoted by  $q, k_I, q_I$ and  $\la$,  for the standard two dimensional Brownian motion $B(t)$. In doing  this it is recalled that the 
  interval  $I(n)$ is defined to be $(-n+1/2, n-1/2)$ instead of $[-n+1,n-1]$,   which  makes difference  in the associated probabilities of the Brownian motion. Put $L_{\pm}=\{t\in \R: \pm t>0\}$ and $\tau_{L_\pm}=\inf\{t>0:  B(t)\in L_\pm\}$, and  define
 $$h^{\pm}_{x}(s)=P_{x}^{BM}[B(\tau_{L_{\pm}})\in ds]/ds \quad \quad (x\in L_{\mp}, \pm s>0),$$
where $P_z^{BM}$ denotes the law of $B(t)$ starting at $z$.
 Then for real  $x> n_*, \,y>- n_*$,
$$q(x,y)=\int_{-\infty}^{-n_*}h^-_{x-n_*}(s-n_*)h^+_{s+n_*}(y+n_*)ds,$$
$$q_{I}(x,y)=q(x,y){\bf 1}(-n_*<y < n_*),$$
$$ k_I(x,y)=h^-_{x-n_*}(y-n_*){\bf 1}(-n_*<y <n_*)$$
and
$q^k$ and $\la$ are given in analogous ways; in particular $q^1=q$ and
$$\la(x,y)=\sum_{k=1}^ \infty q^k(x,y){\bf 1}(y > n_*).$$
 We know that
\beqn\label{*}
h^-_{x}(s)= \frac{\sqrt x}{\pi\, (x-s)}\cdot \frac1{\sqrt{-s}}    \quad (x>0, s<0),
\eeqn
 $$h^{I(n)}_x( s)= \frac{\sqrt{x^2-n_*^2}}{\pi\, (x-s)}\cdot \frac1{\sqrt{n_*^2-s^2}} \quad\quad (x>n_*,  -n_*<s<n_*).$$

The function  $Q$ is extended to
that of reals by 
$$\textstyle Q(u,v) =Q(x,y) \quad \mbox{ for} \quad (u,v)\in (x-\frac12, x+\frac12]\times (y-\frac12, y+\frac12],$$
and similarly for $\La$ and $K_I$. (The summation  in (\ref{QQ}) can  be  then  replaced by  the integration over $y>n_*$.)
With  $Q$ thus extended put $$\eta =Q-q.$$ 

We shall prove the  relations (I) through (VII) given below.  The
symbol $f\asymp g$ means that $C^{-1}g\leq f \leq Cg$. Here and in what follows  $C$ denotes a positive constant which may depend on the law $P_0[S_1=\cdot]$ but  is independent of  any variables $x, n, y, s$ contained therein explicitly or inexplicitly and may change from line to line.  The products of two functions (of two variables) are understood to be that of integral operators in an  analogous way to  (\ref{QQ}): e.g., 
\beqn\label{prod}
\eta q(x,y)=\int_{n_*}^\infty \eta(x,u)q(u,y)du
\eeqn 
(the range of the integration is always the half-line $u>n_*$ where  $q(u, \cdot)$ is defined). 

Let $x>n_*,\,  y>-n_*$ and $-n_*<s<n_*$; $x, y, s$ are  real numbers in (I) through (III).
 
\vskip3mm
(I)  $ \quad \displaystyle q(x,y)\asymp \frac{\sqrt{x-n_*}}{\sqrt{n_*+y}}\cdot \Bigg|\frac1{x-y}\log\frac{x+2n}{y+2n}\Bigg|.$
\vskip3mm\n
The function $t^{-1}\log (1+t)$ is understood to be continuously extended to
$t=0$.  On using  the inequalities   $1/b<(b-a)^{-1}\log (b/a) <1/a$  ($0<a<b$) we infer that  $(x-y)^{-1}\log[(x+2n)/(y+2n)]\asymp x^{-1}$ if
$|x -y|\leq 3n$, which combined with (I)  yields the bound of $q_I$ given in the next item where we also display the explicit form of $k_I$ for convenience.
\vskip3mm
(I$'$) $\quad\displaystyle   q_I(x,s)\asymp \frac{\sqrt{x-n_*}}{\sqrt{n_*+s}}\cdot \frac1{x}\Big(1+\log \frac{x}{\,n_*}\Big);$
\vskip2mm
$\quad\quad\quad \displaystyle k_I(x,s) = \frac1{\pi}\frac{\sqrt{x-n_*}}{(x-s)\sqrt{n_*-s}}.$

\vskip4mm
(II) $\quad\displaystyle \frac {|\eta(x,y)|}{q(x,y)}$
\begin{eqnarray*}
&& = \,  o(1) \quad \mbox{as \quad $(x-n)\wedge(n+y)\to\infty$} \quadd \,\,\mbox{if \,\,\,$\de=0$},\\
&&\leq \, \frac{C}{\sqrt{(x-n_*)\wedge(n_*+y)}} \quadd\quadd \qquad  \mbox{if \,\, $\de>\frac12$.} \,\,\,
 \end{eqnarray*}
 
\vskip4mm
(III) $\quad  |\eta|({\bf 1}+ \la)(q_I+k_I)(x,s)$
\begin{eqnarray*}
 && = \,\Bigg[\frac{1}{\sqrt{n_*^2-s^2}}\wedge h^{I(n)}_x(s)\Bigg]\times o(1)
\quad \mbox{as \quad   $n\to\infty$} \quad  \mbox{if \,\,\,$\de=0$},\\
&& \leq \,  C\frac{1+(\log x/n_*)^2}{\sqrt{x}\sqrt{n_*^2-s^2}} \,\,\quad\quadd  \quadd  \quad  \quad  \quad     \mbox{if\,\, \, $\de>\frac12$.}
 \end{eqnarray*}

\vskip1mm\n
Here $|\eta|$  in (III) stands for the integral  operator of   kernel  $|\eta(x,y)|$ as in (\ref{prod}), ${\bf 1}$,  also  in (III),   for  the identity operator,  and $\de$ in (II) and (III) for the constant  in (\ref{de}). 

 \vskip3mm 
For integers $x\geq n, -n<s<n$,
 \vskip3mm
(IV) $\quad \displaystyle  \sum_{y=n}^\infty \La(x,y) \leq C\sqrt{\frac{x-n_*}{x}},$
 \vskip3mm 
(V) $\,\quad \displaystyle \sum_{y=n}^\infty \La (x,y)\frac1{y-s}
 \leq  C \frac{1}{ n}\sqrt{\frac{x-n_*}{x}}\cdot\log\frac{3n}{n-s},$ 
 \vskip3mm
(VI) $\quad \displaystyle  \sum_{y=n}^{n+N} \La(x,y) \leq C\sqrt{\frac{x-n_*}{x}}\cdot\frac{N}{n}   \quad \quad (N=1,2,\ldots).$
 \vskip3mm
There exist   functions  $\e_j(t)$, $j=1,2$,  of (a single variable) 
such that as $t\to\infty$, $\e_j(t)= o(1)$ or $O(1/\sqrt {t}\,)$  according as $\de=0$  or $\de>1/2$ in (\ref{de}) and that
\vskip3mm
(VII) $ \quad \displaystyle
|k_I-K_I|(x,s)\leq  k_I(x,s)[\e_1(x-n_*)+\e_2(n-s)].$

\vskip5mm 
 The proofs of these results are  postponed to the next section. In the rest of this section  we prove Theorem \ref{thm1} taking them for granted.  
 
By symmetry we may suppose $x\geq n$.  From the identity $$H_x^{I(n)}=({\bf 1}+\La)(K_I+Q_I)(x,\cdot)$$
 and a similar one for $h_x^{I(n)}$ it follows that
$$H_x^{I(n)}-h_x^{I(n)}=({\bf 1}+\La)(K_I-k_I+Q_I-q_I)+(\La-\la)(k_I+q_I)(x,\cdot).$$
Writing  $Q=\La-\La Q$ and $q=\la-q\la$ (valid on $[n_*,\infty)^2$) one finds the identity $\La q- Q\la=\La \eta\la$,  which yields
 \beqn\label{La-la}
 \La-\la = \eta + \La \eta +\eta\la + \La\eta \la=({\bf 1}+\La)\eta({\bf 1}+\la).
 \eeqn
Let $q_I+k_I$ act on the both sides from the right. Let $x$ and $s$ be  integers such that $x\geq n$ and $-n<s<n$. Using  (III), (IV)
and the
simple inequality
$$x^{-\a}\log x/n \leq(e\a)^{-1}n^{-\a} \quad \quad (\a>0),$$
 first observe that
$$\La |\eta|({\bf 1} +  \la)(q_I+k_I)(x,s)\leq \sqrt{\frac{x-n_*}{x}}\cdot\frac{\e(n)}{\sqrt{n^2-s^2}},
$$
where  $\e(t)$  is a  function  of the same meaning as $\e_j(t)$ in (VII),
and then, further using (\ref{La-la}) and  (III), that 
\begin{eqnarray}\label{(6)}
|\La-\la|(q_I+k_I)(x,s)&\leq&  \Bigg[\e(n)\sqrt{x-n_*}+C\bigg(1+\log \frac{x}{n}\bigg)^2\,\Bigg]\frac{1}{\sqrt{x(n^2-s^2)}}\nonumber\\
&\leq& C'\Bigg(\e(n)+\frac1{\sqrt{x-n_*}}\Bigg)h^{I(n)}_x(s).
\end{eqnarray}
The last inequality in particular implies
\beqn\label{(6_1)}
\La (k_I+q_I)(x,s)\leq  Ch^{I(n)}_x(s).
\eeqn
 \vskip3mm

Let  $\de>1/2$ in (\ref{de}). Combined with  (V) and (\ref{(6_1)}) the  bound (VII) shows that 
$$\La |k_I-K_I|(x,s)\leq  C\Bigg( \frac{1}{ n}\sqrt{\frac{x-n_*}{x}}\,\log\frac{3n}{n-s}+h^{I(n)}_x(s)\Bigg)\frac1{\sqrt{n-s}},$$
but we have $n^{-1}\log [3n/(n-s)]\leq 1/\sqrt{n(n-s)}$
so that 
\beqn\label{1-3}
\La |k_I-K_I|(x,s)\leq  C\Bigg(\sqrt{ \frac{x-n_*} {nx(x-s)}}\,+ h^{I(n)}_x(s)\Bigg)\frac1{\sqrt{n-s}}\leq C'h^{I(n)}_x(s)\frac1{\sqrt{n-s}}.
\eeqn
On the other hand by (II) and (I$'$) we have
\beqn\label{q-Q}
|q_I-Q_I|(y,s)\leq C\Bigg(\frac1{\sqrt{y-n}}+\frac1{\sqrt{n+s}}\Bigg)\cdot q_I(y,s)\leq  C\frac{1+\log\,y/n}{\sqrt {y}}\cdot\frac1{n+s}
\eeqn
and   $\La |q_I-Q_I|\leq C\sqrt{x-n_*}/\sqrt {xn}\,(n+s)\leq C'h^{I(n)}_x(s)/\sqrt{n+s}$, which in conjunction  with (\ref{1-3}) gives
\beqn\label{(7)}
({\bf 1}+\La)\Big(|q_I-Q_I|+|k_I-K_I|\Big)(x,s)\leq Ch^{I(n)}_x(s)\Bigg( \frac{1}{\sqrt{n-s}}+ \frac{1}{\sqrt{n+s}}\Bigg).
\eeqn
The bounds (VII),
 (\ref{(6)}), 
 (\ref{q-Q}) and (\ref{(7)}) together yield
 the formula of Theorem \ref{thm1} in the case $\de>1/2$.

\section{Proofs of (I) through (VII)}
\vskip2mm\noindent
{\it Proof of} (I). \,\,Let $x\geq n_*$ and $y>-n_*$. It follows from (\ref{*}) that
\begin{eqnarray}
\label{10}
q(x,y) &=&\frac1{\pi^2}\int_{-\infty}^{-n_*}\frac{\sqrt{x-n_*}}{(x-u)\sqrt{-u+n_*}}\frac{\sqrt{-u-n_*}}{(y-u)\sqrt{n_*+y}}du \nonumber\\
&=&\frac1{\pi^2}\sqrt{\frac{x-n_*}{y+n_*}}J_n(x+n_*,y+n_*),
\end{eqnarray}
where 
$$J_n(a,b)=\int_0^\infty\frac{\sqrt t \,dt}{\sqrt{t+2n_*}(t+a)(t+b)} \quad (a=y+n_*>0,\,b=x+n_*>2n_*).$$

If $n_*\leq a < b$, then
\beq
J_n(a,b)&=&\frac1{a}\int_0^\infty\frac{\sqrt t \,dt}{\sqrt{t+2n_*/a}\, (t+1)(t+b/a)}\\
&\asymp& \frac1{b}+\frac1{a}\int_1^\infty\frac{dt}{(t+1)(t+b/a)} \\
&=&\frac1{b}+\frac1{b-a}\log\frac{a+b}{2a} \\
&\asymp& \frac1{b-a}\log\frac{b}{a}.
\eeq
This  shows (I) in the case $y\geq n_*$. 

In the case $-n_*<y<n_*$, a similar computation gives
\beq
J_n(a,b)=\frac1{2n_*}\int_0^\infty\frac{\sqrt t \,dt}{\sqrt{t+1}\, (t+a/2n_*)(t+b/2n_*)}\asymp \frac1{b}+\frac1{b-a}\log\frac{2n_*+b}{2n_*+a}.
\eeq
Thus  (I) has been proved.

\vskip2mm\noindent
{\it Proof of} (II). \,\, First suppose that $\de>1/2$ in (\ref{de}). 
Then it is shown in \cite{U} (Theorem 1.3) that there exists a constant $C$ such that for $x\geq n$ and $s<n$, 
\beqn\label{thmA}
\Bigg|H_{x-n}^-(s-n)-\frac1{\pi}\cdot\frac1{x-s}\sqrt{\frac{ x-n_*}{n-s}}\,\Bigg|\leq Ch^{-}_{x-n_*}(s-n_*)\Bigg(\frac1{\sqrt{n_*-s}}+\frac1{\sqrt{x-n_*}}\Bigg).
\eeqn
In making application of this and its obvious analogue for $H^+$ there arise four terms to be estimated for computation of the difference $Q-q=H^-H^+-h^-h^+=(H^--h^-)H^++h^-(H^+-h^+)$  (the right side of (\ref{thmA}) is counted two terms), which are equal to those obtained by inserting the factors
$$\frac1{\sqrt{x-n_*}}+\frac1{\sqrt{-u+n_*}} \quad \mbox{and} \quad\frac1{\sqrt{n_*+y}}+\frac1{\sqrt{-u-n_*}}$$
under the integral symbol of the integral of (\ref{10}).
  Among them  only  two terms  require computation, which we  are to show to be not larger
than the sum of the other two.  To this end,
we make the same change of variables that led to the second equality of (\ref{10}) and find that it suffices in view of (I) to verify the following inequalities
$$\int_0^\infty\frac{\sqrt t\, dt}{(t+2n_*)(t+a)(t+b)}\leq \int_0^\infty\frac{dt}{\sqrt{t+2n_*}\,(t+a)(t+b)}\leq \frac{\pi}{(a\vee b)\,\sqrt{a\wedge b}}$$
($a=y+n_*>0,\,b=x+n_*$ as before). The first one is trivial. The second one is
verified by letting  $(a\vee b)^{-1}\int_0^\infty [\sqrt t\, (t+a\wedge b)]^{-1}dt$ dominate the  integral  in the middle.  This completes the proof in the case $\de>1/2$.

The case $\de=0$ is similarly dealt with based on the corresponding result for 
$H_x^-(s)$ (Theorem 1 of \cite{U}).
\vskip2mm\noindent
{\it Proof of} (III). Consider the case $\de>1/2$. Set
$$A(y)=\frac1{\sqrt{x-n_*}}q(x,y)h^{I(n)}_y(s)\sqrt{n_*^2-s^2}, $$
$$   B(y)=q(x,y)\frac1{\sqrt{y+n_*}}h^{I(n)}_y(s)\sqrt{n_*^2-s^2}$$
and
$$I_A=\int_{(x-2n)\vee n}^\infty    A(y)dy \quad \mbox{and}
\quad I_B=\int_{n_*}^{(x-2n)\vee n}B(y)dy.$$
Notice that $1/\sqrt{x-n_*}\leq1/\sqrt{y+n_*}$ if and only if $y\leq x-2n_*$
and that  according to (II)
$$|\eta| h^{I(n)}(x,s)\leq C(I_A+I_B)/\sqrt{n_*^2-s^2}.$$

By (I) 
$$A(y)\asymp \frac{\sqrt{y-n_*}}{y-s}\cdot\frac1{x-y}\log\frac{x+2n}{y+2n}.$$
A simple computation shows that
$$\int_{(x-2n)\vee n}^{(x-2n)\vee (2n)}A(y)dy=O(1/\sqrt x)$$
and
$$\int_{(x-2n)\vee (2n)}^\infty A(y)dy\leq C\int_{(x-2n)\vee (2n)}^\infty\frac1{y^{3/2}}\Bigg(1+\log\frac yx\Bigg)dy =O(1/\sqrt x).$$
Thus $I_A=O(1/\sqrt x)$ (uniformly in $s<n_*$).

As for $I_B$ first we observe
$$B(y)\asymp \frac{\sqrt{x-n_*}}{\sqrt{y+n_*}}\frac{\sqrt{y-n_*}}{y-s}\cdot\frac1{x-y}\log\frac{x+2n}{y+2n}.$$
Let $x\geq 3n$.  It is easy to see that 
$ \int_n^{2n}B(y)dy=$ const $(1/\sqrt x)\log (x/n_*),$
while 
\beq
 \int_{2n}^{x}B(y)dy &\leq& C\sqrt x \int_{2n}^x\frac1{y}\cdot\frac1{x-y}\log\frac{x+2n}{y+2n}dy \\
&\leq& C\frac{2\sqrt x}{x+2n}\int^1_{4n/(x+2n)}\frac1{u(1-u)}\log\frac1{u}\,du
\eeq
and   
the last member is dominated by a constant multiple of $(1/\sqrt x)(\log x/n)^2$ owing to the equality $\int_a^1 u^{-1}\log(1/u)du=\frac12 |\log a|^2$ ($0<a<1$).
 This verifies (III) when $\de>1/2$.  For the case $\de=0$ the same argument as above leads to the upper bound $o(1)/\sqrt{n_*^2-s^2}$; the identity $qh^{I}_{\cdot}=h_\cdot^I-(k_I+q_I) \,(\leq h_\cdot^I)$ gives the other bound $h_\cdot^I\times o(1)$. The proof of (III) is complete.

 \vskip2mm\n
 {\it Proof of} (IV).  Put $p_n=\sup_{x, y\geq n}Q(x,y)$. Then $p_n=O(1\sqrt n)$ and 
$$\sum_{y\geq n}Q^k(x,y)\leq p_n^{k-1}\sum_{y\geq n}Q(x,y)\leq p_n^{k-1}\sum_{y\geq -n}H^-_{x-n}(y)\leq  p_n^{k-1}C\sqrt{\frac{x-n_*}{x}},$$
 hence $\sum_{y=n}^\infty \La(x,y)\leq \sum_k\sum_{y\geq n}Q^k(x,y)\leq C'\sqrt{(x-n_*)/x}$.

 \vskip2mm\noindent
\begin{Lem}\label{lemV_1}\,\,Uniformly for  integers $x\geq n, -n<s<n$,
$$\sum_{y=n}^\infty Q(x,y)\frac1{y-s}\asymp \frac{\sqrt{x-n_*}}{x\sqrt n}
\Bigg( 1+\log \frac{x}{ n}\Bigg)\log\frac{3n}{n-s}.$$ 
\end{Lem}

\noindent
{\it Proof}. By employing (I) and the bound $(x-y)^{-1}\log\frac{x+2n}{y+2n} \asymp x^{-1}[1+\log (x/y)]$ valid for $n\leq y\leq x+n$ one sees that
\beq
\sum_{y=n}^{2n} Q(x,y)\frac1{y-s}
&\asymp& \int_n^{2n}\frac{\sqrt{x-n_*}}{(y-s)\sqrt{n+y}\, }\cdot \frac1{x-y}\log\frac{x+2n}{y+2n} dy\\
&\asymp& 
\frac{\sqrt{x-n_*}}{x\sqrt n}
\Bigg( 1+\log \frac{x}{ n}\Bigg)\log\frac{2n-s}{n-s}
\eeq
as well as
$$\sum_{y=2n+1}^{\infty} Q(x,y)\frac1{y-s}
\leq C
\frac{\sqrt{x-n_*}}{x\sqrt n}\Bigg( 1+\log \frac{x}{ n}\Bigg)
$$
(break the summation according as $y$ is larger  than  $x\vee (2n)$ or not  and  consider the cases $x\leq 2n$ and $x>2n$ separately). These 
 together yield the estimate of the lemma. \qed

\noindent
{\it Proof of} (V).  Letting  $2\sqrt{x/n}$  dominate $1+\log (x/n)$  in the right-hand side of the asymptotic formula of Lemma \ref{lemV_1}   and employing (IV), we have 
$$\sum_{y=n}^{\infty}\La(x,y)\frac1{y-s}=\sum_{y=n}^{\infty}(Q+ \La Q)(x,y)\frac1{y-s}\leq  C\sqrt{\frac{x-n_*}{x}}\,\frac{\log[3n/(n-s)]}{n}.$$
Thus (V) is proved.
\vskip2mm\noindent
{\it Proof of} (VI). As in the proof of Lemma \ref{lemV_1} we have 
$$\sum_{y=n}^{n+N}Q(x,y)\leq C\sqrt{(x-n_*)/x}\, N/n;$$ the estimate of (VI) is then follows from (IV) as in the preceding proof.
\v2\n
{\it Proof of} (VII): This follows immediately from (\ref{thmA}) if $\de>1/2$. In  the case $\de=0$ use  (\ref{(003)}).

\vskip2mm\n
{\it Proof of Proposition \ref{prop0}}. Both formulae (\ref{B-1}) and (\ref{B-2}) of Proposition \ref{prop0} are proved in a similar way as Theorem \ref{thm1} and in  their proofs given below   we omit  details. We first show the   deduction of (\ref{B-2}) from (\ref{B-1}). For the symmetric simple random walk the  right-hand side of (\ref{thmA}) can be replaced by $Ch^{-}_x(s)[(-s)^{-1}+(x\vee 1)^{-1}]$ (cf. \cite{U}) and accordingly we  deduce that
$$(|k_I-K_I|+|q_I-Q_I|)(y,s)\leq \frac{C k_I(y,s)}{[(y-n_*)\wedge(n-s)\wedge(n+s)]},$$ 
$$|\eta|(x,y)\leq \frac{Cq(x.y)}{(x-n_*)\wedge(n_*+y)},$$
$$\sum_{y=n}^\infty \frac{Q(x,y)}{(y-s)\sqrt{y-n_*}}\asymp \frac{\sqrt{x-n_*}}{x\sqrt n}\Bigg(1+\log \frac x{n}\Bigg)\frac1{\sqrt{n-s}},$$
$$\sum_{y=n}^\infty \frac{\La(x,y)}{(y-s)\sqrt{y-n_*}}
\leq \frac{C\sqrt{x-n_*}}{n\sqrt {x}\,\sqrt{n-s}}$$
and 
$$
\sum_{y=n}^\infty \frac{\La(x,y)k_I(y,s)}{y-n_*}
\leq \frac{C\sqrt{x-n_*}}{\sqrt {x}\,n(n-s)},$$
and with these bounds we can proceed as above to obtain (\ref{B-2}).

{\it Proof of (\ref{B-1})}:
The proof is based on an asymptotic expansion of the potential function of the random walk  (cf. \cite{FU}, \cite{U0}, \cite{KS}) from which an application of the reflection principle immediately yields
$$H_{im}(s)=\frac{|m|}{\pi(|s|^2+m^2)}+O\Bigg(\frac1{|s|^3+|m|^3}\Bigg) \quad (m\neq 0),$$
where $H_{im}(s)$ stands for the probability that the first visit to the real axis of the simple random walk starting at $im\in i\Z$  takes place at $s\in \Z$. We proceed as in Section 2. Bearing symmetry of the walk in mind,  this time  we define for $y\in \Z$ and $x\geq 0$,
$$Q(x,y)=\sum^{\infty}_{m=-\infty}H_{ix}(m)H_{im}(y) \quad \mbox{if\,\, $x>0$\quad and \quad $=H_0(y)$\quad  if \,\,$x=0$}$$
and  inductively
$
Q^k(x,y)=\sum_{u=0}^\infty Q^{k-1}(x,u)Q(u,y) \,(k=1,2,\ldots)$. 
We have the corresponding quantities $h$, $h^-$, $q$ and $q^k$ for the standard Brownian motion. Then for $s<0$,\, $x\geq 0$,
$$H_x^-(s)=\La(x,s):=\sum_{k=1}^\infty Q^k(x,s), \quad\quad h_x^-(s)=\la(x,s):=\sum_{k=1}^\infty  q^k(x,s).$$ 
We know that $C^{-1}\la\leq \La\leq C\la$ for some constant $C>0$ (cf. \cite{U}).  We suitably extending $Q$ to the real variables and  put $\eta=Q-q$ as before.  An elementary computation then gives  in turn
$$\eta(x,y)=O(|xy^2|_+^{-1}\wedge |x^2y|_+^{-1}),$$
$$|\eta|\la(x,y)=O(|xy|_+^{-3/2}\wedge |x^2y|_+^{-1}),\quad \La|\eta|(x,y)=O(|x^{1/2}y^2|_+^{-1}) $$
and 
$$
\La|\eta|\la(x,s) \leq \frac{C}{|x^{1/2} s^{3/2}|_+} \asymp \frac{1}{x-s}\sqrt {\frac{x\vee1}{-s}}\times\Bigg[\frac{1}{{-s}}+\frac1{{x\vee1}}\Bigg] \quad (s\leq -1,\,x\geq 0),\,$$
where $|a|_+=|a|\vee 1$. Thus  (\ref{B-1}) follows because of  the identity 
$$\La-\la=(1+\La)\eta(1+\la).$$

The proof of Proposition  \ref{prop0} is complete.  \qed

\section{Estimation of $H^{I(n)}_x(s)$ near the edges}

We continue the arguments of the preceding section to estimate $H^{I(n)}_x(s)$ mainly in the case when $\de>1/2$ and either $n-s$ or $n+s$ is small in comparison with $x-n$. The case when $\de=0$ or $x-n$ is not large can be similarly dealt with and is only briefly discussed at the end of this section.   
 
 Theorems \ref{thm5} and \ref{thmII-2} given below are based on  the following result from \cite{U} ((10) of Theorem 1.3 and its dual): if $\de>1/2$ in (\ref{de}), then for $x\geq s>0$, 
 \beqn\label{thmA0}
 H_x^-(-s)=\frac{\sqrt x}{\pi (x+s)}\,\mu^-(s)+O\Big(\frac {1}{x}\Big),
  \eeqn
 \beqn\label{thmA1}
 H_{-x}^+(s)=\frac{\sqrt x}{\pi (x+s)}\,\nu^-(s)+O\Big(\frac {1}{x}\Big),
  \eeqn
where  $\mu^-(s)=\mu(-s)$ and $\nu^-(s)=\nu(-s)$. The following theorem concerns particularly to the case when $(x-n)/(n-s)\to\infty$ so that 
$h^{I(n)}_x(s) >> (x-s)^{-1}$.

\begin{Thm}\label{thm5}\,\,If $\de>1/2$ in (\ref{de}), then uniformly for integers $n>1, 0\leq s<n$ and $x\geq n$,
$$H^{I(n)}_x(s)=\sqrt{n_*-s}\,\mu(-n+s)h^{I(n)}_x(s)+O\Bigg(\frac{\log n}{n} +\frac1{x-s}\Bigg).$$
\end{Thm}

\vskip2mm\noindent
{\it Proof}. Make decomposition $H_x^{I(n)}=K_I+Q_I+\La (K_I+Q_I)$ and infer from (\ref{thmA0}) that
\beqn\label{K+Q}
K_I=\sqrt{n_*-s}\,\mu^-(n-s)k_I+O(1/(x-s)),
\eeqn
$$\La K_I(x,s)=\sum_{y=n}^\infty \La(x,y)\Bigg[\sqrt{n_*-s}\,\mu^-(n-s)k_I(y,s) +O\bigg(\frac1{y-s}\bigg)\Bigg].$$
In view of (I) we have
$\sup_{s\geq 0, x>n_*} q_I(x,s)\leq C/n,$
which in particular shows that  
$$\la q_I(x,s)=O(1/n) \quad \mbox{ uniformly for} \quad 0\leq s<n_*, x>n_*.$$ 
Thus, on employing (I$'$), for $s>0$, 
\beqn\label{19}
q_I+\la q_I =O(1/n) \quad \mbox{and} \quad Q_I+\La Q_I\leq C(q_I+\la q_I)=O(1/n).
\eeqn

By (V) of the preceding section we have
\beqn\label{2-1}
\sum_{y=n}^\infty \La(x,y)\frac1{y-s} \leq C\sqrt{\frac{x-n_*}{x}}\,\frac{\log n}{n},
\eeqn
so that
$$
\La K_I(x,s)=\sqrt{n_*-s}\,\mu^-(n-s)\La k_I(x,s)+O(n^{-1}{\log n} ).$$
Here the factor $\sqrt{(x-n_*)/x}$ on the right side of (\ref{2-1}) is replaced by 1: the loss of accuracy to  the estimate of $H_x^{I(n)}$  caused by this replacement is  small in comparison with the  error term $O(1/(x-s))$ in (\ref{K+Q}).
By (\ref{(6)})
$
\La k_I=\la k_I+(\La-\la)k_I
=\la k_I+O(1/n);$
 hence 
$$\La K_I=\sqrt{n_*-s}\,\mu^-(n-s)\la k_I+O(n^{-1}\log n),$$
which together with (\ref{K+Q}),  (\ref{19})  yields the assertion of the theorem.
  \qed
 \begin{Thm}\label{thmII-2}\,\,If $\de>1/2$ in (\ref{de}), then uniformly for integers $n>1, -n< s<0$ and $x\geq n$, 
\beq
H^{I(n)}_x(s)&=&\sqrt{n_*+s}\,\nu(-n-s)h^{I(n)}_x(s)\\
&& \times\,\left[ 1+O\Bigg(\sqrt {\frac{s+n_*}{n}}\cdot\log n\Bigg)  +O\Bigg(\sqrt{\frac{x}{n(x-n_*)}}\, \Bigg)\right].
\eeq
\end{Thm}
\vskip2mm\noindent
{\it Proof}.  We  make decomposition
\begin{eqnarray}\label{eq-thm8}
H^{I(n)}_x(s)&=&K_I(x,s)+\sum_{y=-\infty}^{-n}H^-_{x-n}(y-n)H^{I(n)}_y(s)  \nonumber \\
&=&K_I(x,s)+\sum_{y=-\infty}^{-n}(H^-_{x-n}-h^-_{x-n})(y-n)H^{I(n)}_y(s) \nonumber \\
&&\,\, + \sum_{y=-\infty}^{-n} h^-_{x-n}(y-n)H^{I(n)}_y(s). 
\end{eqnarray}
For evaluation of the second  sum of the last line we substitute the estimate of Theorem \ref{thm5} for $H_x^{I(n)}(s)$ (with $S_n$ replaced by $-S_n$, hence $\sqrt{n_*-s}\, \mu^-(n-s)$ by $\sqrt{n_*+s}\, \nu^-(n+s)$), use the expression of  $h^{I(n)}_x(s)$ analogous to the first expression in (\ref{eq-thm8}) of $H^{I(n)}_x(s)$ and observe  that 
$$\sup_{-n<s\leq 0}k_I(x,s)\asymp\frac{\sqrt {x-n_*}}{x\sqrt{n}};\,\,\, \int_{-\infty}^{-n} h^-_{x-n}(y-n)dy\asymp\frac{\sqrt {x-n_*}}{\sqrt x};$$
$$\quad \mbox{and} \quad \int_{-\infty}^{-n} h^-_{x-n}(y-n)(s-y)^{-1}dy\asymp \frac{\sqrt {x-n_*}}{x\sqrt{n}}\cdot \log \frac{n}{s+n} \quad$$
  to obtain
\beqn\label{EST1}
\sum_{y=-\infty}^{-n} h^-_{x-n}(y-n)H^{I(n)}_y(s)=\sqrt{n_*+s}\,\nu^-(n+s)h^{I(n)}_x(s)+O\Bigg(\frac{\sqrt {x-n_*}}{\sqrt{x}\,n}\cdot \log n\Bigg).
\eeqn
For evaluation of the first sum   apply   (\ref{thmA}) to have    
$$|H^-_{x-n}(y')-h^-_{x-n}(y')|\leq C((x-n_*)^{-1/2}+n^{-1/2})h^-_{x-n}(y')$$
 for $y'\leq -2n$; also use the bound $H^{I(n)}_x(s)\leq  Ch^{I(n)}_x(s)$ that follows from Theorem \ref{thm1}. These bounds as well as $\int_{-\infty}^{-n_*}h^-_{x-n}(y-n)h^{I(n)}_y(s)dy\leq h^{I(n)}_x(s)$ yield 
\beq
\sum_{y=-\infty}^{-n}|(H^-_{x-n}-h^-_{x-n})(y-n)|H^{I(n)}_y(s) &\leq&  C\Bigg(\frac{1}{\sqrt{x-n_*}}+\frac1{\sqrt n}\Bigg)h^{I(n)}_x(s) \\&\asymp&  \frac{\sqrt x}{\sqrt{(x-n_*)n}}h^{I(n)}_x(s).
\eeq
Combined with (\ref{eq-thm8}) and (\ref{EST1}) this completes the proof of the theorem. \qed 

As being mentioned at the beginning of this section our estimation of $H^{I(n)}$ made above is appropriate if $x-n$ is large in comparison with $n\pm s$. When $x-n$ is not large, it is better to replace $h^-_{x-n}(y-n)$ by
$$\frac{\sigma^2 \nu(x-n)}{2\sqrt{x-n_*}}\,h^-_{x-n_*}(y-n)$$
in  (\ref{eq-thm8}); also make use of the corresponding estimate of $H^{\pm}_{x}$ in \cite{U}. 

\vskip2mm\n
{\it Proof of Theorem \ref{thm2}} ({\it the case $\de=0$}). The proof is based on an estimate of $H^-_x(s)$ verified in \cite{U} (Theorem 1.1). The case when $(x-n)\wedge (n-s)\wedge (n+s)\to\infty$ the assertion is included in Theorem \ref{thm1}. The other case  is dealt with as in the proofs of 
 Theorems \ref{thm5} and \ref{thmII-2} by employing (VI).  The computations to be carried out this time are
much simpler in either case.

\section{Proof of Theorem \ref{thm21}}
We always have the relation
\beqn\label{5.0}
H^{I(n)}_x(s) =H_x(s) + \sum_{|x_1|\geq n} H_x(x_1)H^{I(n)}_{x_1}(s)   \quad\quad s\in I(n),
\eeqn
where  $H_x(s)$, $x, s\in \Z$,  denote  the probability that $s$ is the  site where the real axis is  hit  for the first time after time 0 by the walk $S_n$ started at $x$. Suppose $E_0[\,|Y|^2\log |Y|\,]<\infty$.  Then
\beqn\label{5.00}
H_x(s)=H_0(s-x)= \frac{\sigma^2}{\pi (s-x)^2}(1+o(1))  \quad \quad \mbox{as} \,\, |s-x| \to\infty
\eeqn
(cf. \cite{U1}).
The proof of  Theorem \ref{thm21} is relatively simple owing to  (\ref{5.00}). Before proceeding with the proof  we mention a few  points to be recognized.  For the estimate of $H^{I(n)}_x(s),  s, x\in  I(n)$  we may suppose that
 $x\geq0$ for obvious reason and then,  $s<x$,  by considering the time-reversed walk. Also  (ii$'$) is a dual statement of (ii) as already noted. There is some possibility of improving the estimates in certain cases that we do not take up in this paper, and for that purpose  some details given below would be helpful .

From what is noticed above we  may suppose 
\beqn\label{5.1}
 x\geq0  \quad \mbox{and} \quad  x-s \to\infty.
\eeqn 
Then, according to  Theorem \ref{thm2}, for $x_1\geq n$, 
\begin{eqnarray}\label{5.21}
H^{I(n)}_{x_1}(s)& =&  \Bigg[ \frac{\nu(-n-s)}{\sqrt{n-s}}\Bigg]  \frac{\sigma^2 \nu(x_1-n) \sqrt{n +x_1}}{2\pi  |x_1-s| }(1+o(1))\\
&=& \Big [ \sqrt{n+s}\,\nu(-n-s) \Big] h^{I(n)}_{x_1}(s)\frac{\sigma^2 \nu(x_1-n)}{2 \sqrt{x_1-n}}(1+o(1)).
\label{5.22}
 \end{eqnarray}
 (This is obtained  directly  from Theorem \ref{thm2} (ii), which also covers  the case $s>0$ under  (\ref{5.1}).)

\n
{\it Proof of {\rm (ii)}}. \,  The Wiener-Hopf equation for $\nu$ may be written as
\beqn \label{5.3} 
\sum_{x_1=n}^\infty H_x(x_1)\nu(x_1-n) =\nu(-n+x)
\eeqn
(see (\ref{WH})). We claim that as $(n-x)/(n-s)\to 0$,
 \beqn \label{5.4}
J_{n,s,x}:= \sum_{x_1\geq n} H_x(x_1)H^{I(n)}_{x_1}(s)  =   \Bigg[ \frac{\nu(-n-s)}{\sqrt{n-s}}\Bigg] \cdot \frac{  \sigma^2 \nu(-n+x) \sqrt {2n} }{2 \,\pi (n-s)} \Big(1+o(1)\Big).
\eeqn 
Put $\xi=n-x$ and observe that owing to (\ref{mu0})   the summation in (\ref{5.3}) may be restricted to $x_1\leq n+K\xi$ by choosing  $K$  large enough .  Then,  on  looking at (\ref{5.21}),  the claim  (\ref{5.4}) follows if we show that for each  $\e>0$ we can find $K>1$ such that
$$\sum_{x_1 \geq n+K\xi } H_x(x_1)H^{I(n)}_{x_1}(s) \leq \e J_{n,s,x}.$$
However,  by simple consideration this  reduces to 
$$\int_{K\xi}^\infty \frac{\sqrt{2n+u}}{(u+\xi)^2(u+n-s)\sqrt u}du \leq  \frac{\e \sqrt n}{(n-s)\sqrt \xi}, $$
which is certainly true if $K$ is large enough. Thus the claim is verified. One can easily check that 
$$ \sum_{x_1\leq -n} H_x(x_1)H^{I(n)}_{x_1}(s) = O\Bigg(\frac1{[(n+s)\wedge \sqrt n\,] n^2}\Bigg) = o(J_{n,s,x}).$$
 Finally notice that $(n-s)/( x-s)\to 1$ and hence the right side of (\ref{5.4}) may be identical to that of the required formula. The proof of (ii) is complete. \qed

 \v2\n
{\it Proof of  {\rm (i)}}.  We may suppose that $(n-x)/(n-s)$ is bounded away from zero, the case $(n-x)/(n-s)\to 0$  being  included in (ii) that is proved above.    Under this condition we see
$$\sum_{n\leq |x_1|\leq n+K} H_x(x_1)H^{I(n)}_{x_1}(s)  =o\bigg(\frac{1}{(s-x)^2}\bigg)$$
 (indeed this is valid if $(n-x)^4/(n-s)\to \infty$; the contribution of $-n - K\leq x_1\leq -n$ is easy to estimate), namely the sum above is negligible if compared with $H_x(s)$.
Hence one can  replace the ratio appearing last in (\ref{5.22}) by 1. Also $\nu(-n-s)$ may be replaced by $1/\sqrt{n+s}$  and,  substituting the resulting expression into (\ref{5.0}) and applying Lemma 
 \ref{lemA1} of Appendix (B), we conclude the formula of (i).  \qed
 
{\sc Remark 2.}  We could have employed the identity
\beqn\label{5.5}
H^{I(n)}_x(s) =H^-_{x-n}(s-n) + \sum_{x_1\geq -n} H^-_{x-n}(x_1-n)H^{I(n)}_{x_1}(s)   \quad\quad s\in I(n),
\eeqn
instead of (\ref{5.0}).  This way is   simpler owing to (\ref{(004)}) except for a tedious computation for exact evaluation  of the definite integral of a  certain rational function.
 
\section{Appendices}

{\bf (A)} \,Let  $D$ be the complement of the line segment with edges  at $\pm1$:
$$D=\C\setminus \{s: -1\leq s\leq 1\}$$
and  denote the Poisson  kernel  (density of harmonic measure)  for $D$ by $h_D(z, s\pm i0)$.
 Putting $h^{[-1,1]}_z(s)= h_D(z, s + i0)+h_D(z, s - i0)$, we have
$$h^{I(n)}_z(s) = h^{[-1,1]}_{z/n_*} (s/n_*)/n_* .$$

We compute  $h_D(z, s\pm i0)$ by using the conformal invariance of harmonic measures. 
The function $z =\frac12 (w+w^{-1})$  univalently maps the exterior of the unit circle 
onto $D$. Denote by $f(z)$ its inverse map, which may be represented
by 
$$f(z)=z+\sqrt {z^2-1},\quad z\in D$$
  with the standard choice of a branch of the square root (so that  $f(\pm s)=\pm s\pm\sqrt {s^2-1}$ for $s>1$  and  $f(s\pm i0)=s \pm   i\sqrt{1-s^2}$  for $-1<s<1$). As  $w=f(z)$ moves on  a circle centered at the origin  counter-clockwise   starting at a point $R>1$,  $z$  describes  the ellipse $[2x/(R+R^{-1})]^2+[2y/(R-R^{-1})]^2=1$ (which  surrounds the segment  $-1\leq s\leq 1$ and shrinks to it as $R\downarrow 1$)  rotating also counter-clockwise and  starting at the point $f(R) =\frac12(R+R^{-1}) \in (1, R)$.  (Cf. \cite{Ah}:p.94 or \cite{N}:p.269).  The Poisson kernel for the exterior of the unit circle is given by 
$$K(Re^{i\th}, e^{i\th'})=\frac{R^2-1}{2\pi (R^2-2 R \cos (\th-\th')+1)},\quad R>1.$$
Put
$$\th(z)=\arg f(z).$$
 Then,  for  $-1\leq s\leq 1$,  $\th(s \pm i0)=\pm \arccos s\in (-\pi,\pi)$, so that $|d\th(s\pm i0)|=ds/\sqrt{1-s^2}$; 
 thus   the conformal  invariance   shows  that  
 \begin{eqnarray}\label{A1}
&&  h_D(z, s\pm i0)  \nonumber\\
&&\,\, =  \frac1{2\pi}\cdot \frac{|f(z)|^2-1}{|f(z)|^2-2|f(z)|\cos \,[\, \th(z)-\th(s\pm i0) ]+1}\cdot \frac1{\sqrt{1-s^2}} \quad \quad \\
  &&\,\, = \frac1{2\pi}\cdot \frac{|f(z)|^2 -1}{|f(z)|^2-2[s \,\Re f(z) \pm \sqrt{1-s^2}\, \Im f(z)]+1}\cdot \frac1{\sqrt{1-s^2}},
 \nonumber \label{A1'}
  \end{eqnarray}
which, for $z=x \in \R\setminus [-1,1]$,  reduces to
 $$h_x^{[-1,1]}(s) =2h_D(x, s+i0) = \frac{\sqrt{x^2-1}}{\pi\, |x-s|}\cdot \frac1{\sqrt{1-s^2}}.$$

 Let $Q=(\sigma_{ij})$ be a $2\times 2$   matrix that is symmetric and positive definite and $\tilde h_D(z, s \pm i0)$  the corresponding hitting density for the process
$Q^{1/2} B_t$, a two-dimensional  Brownian motion of  mean zero  and  the covariance matrix $tQ$.
 Then for $z\in D$,
 \beqn\label{A2}
 \tilde h_D(z, s\pm i0) =  h_{D}(\tilde z, s\pm i0), \quad \tilde z =(x-\om y) +i \la y,
 \eeqn
 where $\om =\sigma_{12}/\sigma_{22}$ and   $\la =\sigma^2/\sigma_{22} =\sqrt{\sigma_{11}/\sigma_{22}-\om^2}$.
If $z$ is real, the identity above  follows  immediately from the rotation invariance of the standard Brownian motion.      In view of the strong Markov property of $S$ its full validity  is deduced  from   the identity  thus restricted in conjunction of the corresponding identity for the Poisson kernel of the upper half plane  (see  (E) below).

\v2\n
{\bf (B)}  \, In Section 5 (at the end of it)  we have used the following lemma.
\begin{Lem}\label{lemA1}\,   For $x, s \in (-n_*, n_*)$ with $s\neq x$, 
\beqn\label{A3}
\frac1{(s-x)^2}+ \int_{|\xi|\geq n_*} \frac1{(\xi-x)^2} h^{I(n)}_\xi(s)d\xi 
= \frac{n_*^2 -xs}{ (s-x)^2\sqrt{(n_*^2-x^2)(n_*^2-s^2)}}.
\eeqn
\end{Lem}
\v2\n
\pf\,  By the scaling property we may suppose the interval  $I(n)$ to be $[-1,1]$.  Let $h_D$ be as in (A) and $h_z(s)$ be the Poisson kernel on the upper half plane: $h_z(s) = y/\pi (y^2+(s-x)^2)$. Then
$$h_z^{[-1,1]}(s) =h_D(z,s + i0)+h_D(z,s - i0) =  h_z(s) + \int_{|\xi|>1} h_z(\xi)h_\xi^{[-1,1]}(s) d\xi,$$ 
which shows  $ \lim_{y\downarrow 0}\pi y^{-1} h_{x+iy}^{[-1,1]}(s) $ equals L.H.S. of (\ref{A3}). The lemma  therefore follows if we verify that if $|x|<1$ and $ |s|<1$,
\beqn\label{A4}
 \lim_{y\downarrow 0} \frac {\pi [h_D(x+iy, s+i0) +\pi h_D(x+iy, s-i0)] }{y} = \mbox{ R.H.S. of  (\ref{A3})}.
 \eeqn
  If $w=- (1-x^2 +y^2) +i 2xy$ and $\phi = \pi - \arg w \in (-\pi/2,\pi/2)$, then
$$|f(z)|^2 =(x+|w| ^{1/2}\sin \phi)^2 +(y+|w| ^{1/2} \cos \phi)^2 $$
and we see that $y^{-1}(|f(z)|^2 -1) \to 2/\sqrt{1- x^2}$.  In view of (\ref{A1}) this shows that
$$ \lim_{y\downarrow 0} \frac  {\pi h_D(x+iy, s\pm i0) }{y} = \frac{1}{2(1-\cos (\th_x\mp \th_s))\sqrt{(1-x^2)(1-s^2)}}$$
where  $\cos \th_t=t$ with $\th_t \in (0,\pi)$  for $-1<t<1$. Now (\ref{A4}) follows from the identity
$$\frac1{1-\cos (\th_x - \th_s)}+ \frac1{1-\cos (\th_x + \th_s)}  = \frac{2-2 \cos \th_x \cos \th_s}{(\cos \th_x - \cos \th_s)^2} =\frac{2(1-xs)}{(x-s)^2}.$$
\qed

\n
{\bf (C)} \,  Let $(X_n)$ be the imbedded  walk on the real axis mentioned in Section 1.  In other words  $(X_n)$ is the one-dimensional random walk with the transition probability $p^X(x,y)= H_0(y-x)$, where $H_0(x)$ be the hitting distribution of the real line for our random walk starting at the origin as being introduced in (\ref{5.0}).  Let
  $ \mu(x)$, $x\geq 0$  be a   renewal  function for the ascending ladder-height variables of the walk $(X_n)$ and $\nu$  its dual;  they are  normalized so as to satisfy  $\lim_{x\to\infty}\mu(x)/\nu(x)=1$ and given by
$$\mu(x)=\frac{\sqrt{\pi}e^{\th_+}}{\sigma}(v_0+\cdots+v_x),\quad
\nu(x)=\frac{\sqrt{\pi}e^{-\th_+}}{\sigma}(u_0+\cdots+u_x)$$
for $x=0, 1, 2, \ldots$  (\cite{S}:p.\,212), where if  $c=\exp\Big(-\sum_{k=1}^\infty\frac1{k}P_0[X_k=0]\Big)$,
\beqn\label{a}
\th_+=\sum_{k=1}^\infty\frac1{k}\Big(\frac12-P_0[X_k>0]\Big)+\frac12 \log c\,=\, \frac12 \sum_{k=1}^\infty\frac1{k}\Big(P_0[X_k<0]-P_0[X_k>0]\Big),
\eeqn
$v_0=u_0=1/\sqrt c$ and 
$\sqrt c v_k$ (resp. $\sqrt c u_k$) equals the probability that  the ascending (resp. descending) ladder-height process visits $k$ (resp.  $-k$) (\cite{S}:p.202, 203).  Then $\mu$ and  $\nu$ 
 are  positive solutions of the Wiener-Hopf integral equations associated with the kernels  $p^X(x,y)$ and $p^X(y, x)$ ($x, y\geq 0$), respectively  to the effect that for $x= 0, 1, 2,\ldots$,
\beqn\label{WH}
\mu(x)=\sum_{k=0}^\infty \mu(k) H_0(x-k)
\quad\mbox{and} \quad
\nu(x)= \sum_{k=0}^\infty \nu(k)H_0(-x+k)
\eeqn
(see \cite{S}:p.332 for uniqueness of positive solutions). 
\v2
\n
{\bf (D)}  \,  Suppose that $E_0[\,|Y|^2\log |Y|\,]<\infty$ (as in Theorem \ref{thm21}). Put  $\sigma_j^2= E[(S_1^{(j)})^2]$ ($j=1,2$) and $\sigma_{12}= E[S_1^{(1)}S_1^{(2)}]$. Then putting for  $s, x,y\in \Z$
$$H_z(s) = P_z[\exists j, S_n=s \,\, \mbox{and}\,\, S_k\notin \R\,\,\mbox{for}\,\, 1\leq k<j] \quad\quad (z=x+iy)
$$
it is shown in \cite{U1} (Theorem 1.2) that  as  $|x-s|+|y|\to\infty$
\beqn\label{A_D1}
H_z(s) =  \frac{\, 1\vee \sigma^2_{2}a(y)}{\pi \|s- z\|^2}(1+o(1)).
\eeqn
(Note that $\sigma_2^2 a(y) \geq |y|$ \cite{S}:P28.8, P31.1.)  Here 
$$\|z\|^2 = \sigma^{-2} (\sigma_{2}^{2} x^2- 2\sigma_{12}xy+\sigma_{1}^{2}y^2)$$
and $a(y), y\neq 0$ is  the potential function of the one-dimensional walk  $S^{(2)}_n$.
Denote  the Poisson kernel on the upper  half plane by $h_z(s)$ as in (B).  Putting  $\la = \sigma^2/\sigma_2^2$ and $\om = \sigma_{12}/\sigma_2^2$ we have $\|z\|^2 = [(x+\om y)^2 + (\la y)^2]/\la$, and for $y\neq 0$, (\ref{A_D1}) is written as 
\beqn\label{A_D2}
H_z(s) =\frac{1\vee \sigma_2^2a(y)}{ |y|} h_{x-\om y +i\la y}(s)(1+o(1)).
\eeqn
Using this   we can readily deduce from Theorems \ref{thm2} and  $2'$  that  as $ |z|\wedge |z-s| \to\infty$
\beqn\label{A_D3}
H^{I(n)}_{z}(s) =\frac{ \sigma_2^2a(y)}{|y|} \mu(-n+s)\nu(-n+s) \sqrt{n_*^2 -s}\, h^{I(n)}_{x-\om y+i\la y}(s)(1+o(1))
 \eeqn
for $y\neq 0$.  

\v2
\n
{\bf (E)}  \, 
   In view of Donsker's invariance principle the relation (\ref{A_D2}) (resp. (\ref{A_D3})) incidentally shows that $h_{x-\om y +i\la y}(s)$ (resp.  $h^{[-1,1]}_{x-\om y +i\la y}(s)$) is  the  density of the hitting distribution of the real line  (resp. the interval $I(n)$) for the process  $Q^{1/2}B_t$.

We give a direct algebraic verification.   We may replace  $B_t$ by  $UB_t$ with any orthogonal matrix $U$ and choose $U$ so that  the matrix $Q^{1/2}U$ 
sends the real line to itself.
A simple algebraic manipulation leads to 
 $$ (Q^{1/2}U)^{-1} = \frac{\sigma_2}{\sigma^2} \left( \begin{array}{cc} 1 & -\om \\ 0 & \la \end{array} \right).$$
 Let  $z=x+iy$, $\tilde z = x-\om y+ i\la y$ and $c=\sigma_2/ \sigma^2$. Noting that $c\tilde z$ corresponds to $z$, we then find, with the notation analogous to $\tilde h_D(z,s)$ in (\ref{A2}),  that 
  $$ \tilde h_{z}( s)=  c h_{c\tilde z}( cs),$$
of which  the right-hand side equals $h_{\tilde z}( s)$,  yielding  the analogue of (\ref{A2}) as desired.

\vskip20mm


\begin{thebibliography}{99}
\baselineskip=9pt

\bibitem{Ah}  L. V. Ahlfors, Complex analysis, 3rd ed.,  McGraw-Hill, 1979 

 \bibitem{F} Y. Fukai, Hitting distribution to a quadrant of two-dimensional random walk, Kodai Math. Jour. {\bf 23} (2000) 35-80.

\bibitem {FU} Y. Fukai and K. Uchiyama, Potential kernel for  two-dimensional random walk, Ann. Probab. {\bf 24} (1992) 1979-1992. 


\bibitem {K} H. Kesten, Hitting probabilities of random walks on $\Z^d$, Stoch. Proc. Appl.
 {\bf 25} (1987) 165-184.
 
\bibitem {K1} H. Kesten, How long are the arms of DLA, J. Phys. A {\bf 20} (1987) L29-L33.

\bibitem {K2} H. Kesten, Some caricatures of the multiple contact diffusion-limitted aggregation and $\eta$-model, pp.179-227 in Stochastic Analysis, eds. M.T. Barlow and N.H. Bingham, Cambridge Univ. Press, 1991.

\bibitem {KS}  G. Kozma and E. Schreiber, An asymptotic expansion for the discrete harmonic potential, Electronic J.P. {\bf 9} (2004) 1-17.

\bibitem {L} G.F. Lawler, Intersections of random walks, Birkh\"auser, 1991.

\bibitem {LL} G. F. Lawler and V. Limic,  Random walks: A modern introduction, Cambridge Univ. Press, 2010.



\bibitem {N} Z. Nehari, Conformal mapping,  McGraw-Hill,  1952 

\bibitem {S} F. Spitzer, Principles of Random Walks, Van Nostrand, Princeton, 1964. 

\bibitem{U0}  K. Uchiyama, Green's functions for random walks on $\Z^d$, Proc. London Math. Soc. {\bf 77} (1998), 215-240.

\bibitem {U}   K. Uchiyama, The hitting distribution of the  half real line for two dimensional random walks. Arkiv f\"or Matematik  {\bf 48}  (2010), 371-393.

\bibitem {Uer}   K. Uchiyama, Erratum to: The hitting distribution of the  half real line for two dimensional random walks. Arkiv f\"or Matematik  {\bf 50}  (2012),  199-200. 

\bibitem {U1}   K. Uchiyama, The hitting distribution of a line for two dimensional random walks, Trans. Amer. Math. Soc. vol. 362, (2010) 2559-2588.  

\end{thebibliography}
\end{document}